\documentclass[12pt]{article}
\usepackage{amsfonts}
\usepackage{color}
\begin{document}

\newtheorem{theorem}{Theorem}[section]
\newtheorem{prop}[theorem]{Proposition}
\newtheorem{cor}[theorem]{Corollary}
\newtheorem{lemma}[theorem]{Lemma}
\newtheorem{ex}[theorem]{Example}
\newtheorem{no}[theorem]{Note}
\newtheorem{unnumber}{}
\renewcommand{\theunnumber}{\relax}
\newtheorem{prepf}[unnumber]{Proof}
\newenvironment{pf}{\prepf\rm}{\endprepf}
\newcommand{\qed}{\qquad$\square$}

\title{Counting false entries in truth tables of bracketed formulae connected by implication}
\author{Peter J. Cameron and Volkan Yildiz\\\\
School of Mathematical Sciences\\
Queen Mary, University of London\\
Mile End Road\\
London E1 4NS\\
\texttt{p.j.cameron@qmul.ac.uk}\\\\
Department of Mathematics\\
King's College London\\
Strand, London WC2R 2LS\\
\texttt{volkan.yildiz@kcl.ac.uk}}
\date{July 14, 2010}
\maketitle

\begin{abstract}
In this paper we count the number of rows $f_n$ with the value ``false'' in
the truth tables of all bracketed formulae with $n$ distinct variables connected by
the binary connective of implication. We find a recurrence and an asymptotic
formulae for $f_n$. We also show that the ratio of $f_n$ to the total number 
of rows converges to $(3-\sqrt{3})/6$.
\end{abstract}

{\footnotesize

}

\pagebreak

\section{Introduction}

In this paper we study enumerative and asymptotic questions on formulae of
propositional calculus which are correctly bracketed chains of implications.

For brevity, we represent truth values of propositional variables and formulae
by $1$ for ``true'' and $0$ for ``false''. 

We begin by stating some important notions of  propositional logic.
The propositional language consists of propositional variables 
$p_1,p_2,\ldots ,p_n$ and  symbols called connectives. 
The well known connectives are `not', `and', `or', `implies', and `if and only if',
which we write as  $\neg$, $\wedge$, $\vee$, $\rightarrow$, and $\leftrightarrow$,  respectively. 
The formulae of propositional logic, are expressions that can be obtained recursively 
from propositional variables by applying connectives. 
More precisely:
\begin{itemize}
\item[(1)] A propositional variable is a formula.
\item[(2)] If $\phi$ and $\psi$ are formulae, then so are $\neg\phi,\; \phi\to \psi,\; \phi\leftrightarrow \psi, \; \phi\wedge \psi,\; \phi\vee \psi$.
\end{itemize}
For unambiguity, brackets are also used in formulae.
For example, we need to be able to distinguish $p_1\to (p_2 \to p_3)$ from $(p_1\to p_2)\to p_3$.
Note that, in U.K., left and right brackets are denoted by the symbols `(' and `)', respectively, 
whereas in U.S., they are denoted by the symbols `[' and `]'. 

Any formula, $\phi$, which involves the propositional variables $p_1,\ldots, p_n$
can be used to define a function of $n$ variables, called `a truth function' 
or `a propositional function', that is, a function from $\{0,1\}^n$ to $\{0,1\}$. 
Since, $|\{0,1\}^n| =|\{0,1\}|^n= 2^n$, the $n$-ary Cartesian product 
$\{0,1\}^n$ has $2^n$ elements.
 Which is the number of rows of a truth table with $n$ 
variables. As is well known, there are $2^{2^n}$ 
propositional functions, each of which can be represented by a formula 
involving the connectives $\neg$, $\vee$ and $\wedge$.

The function represented by a formula is conveniently calculated using a truth
table. Where each row of the truth table corresponds to a valuation.
A valuation is a function $\nu$ from the set of propositions $\{p_1,\ldots,p_n\}$  
to the set $\{1,0\}$. Thus a valuation is an assignment of values to the variables 
$p_1,\ldots ,p_n $, with consequent assignment of values to formulae.

For more information on standard propositional logic the reader can refer to the following books, 
\cite{PJCslc} and \cite{DM}.

We are interested in \emph{bracketed implications}, which are formulae
obtained from $p_1\to p_2\to\cdots\to p_n$ by inserting brackets so that
the result is well-formed, where $p_1,\ldots,p_n$ are distinct propositions. 

The binary connective $\to$  ``implies'' is defined as usual by the rule that, 
for any valuation $\nu$,
\[\nu(\phi\to\psi)=\cases{0 & if $\nu(\phi)=1$ and $\nu(\psi)=0$,\cr\cr
1 & otherwise.\cr}\]
\begin{ex} 
Here are the truth tables, (merged into one), for the two bracketed implications
in $n=3$ variables. Where the corresponding rows with the value false are in blue:

\[
\begin{array}{|l|l|l|c|c|}
\hline p_1 & p_2 & p_3 & p_1\to (p_2 \to p_3) & (p_1 \to p_2)\to p_3 \\
\hline 1 & 1 & 1 & 1 & 1 \\
\hline 1 & 1 & 0 & {\color{blue} 0} & {\color{blue} 0} \\
\hline 1 & 0 & 1 & 1 & 1 \\
\hline 1 & 0 & 0 & 1 & 1 \\
\hline 0 & 1 & 1 & 1 & 1 \\
\hline 0 & 1 & 0 & 1 & {\color{blue} 0 }\\
\hline 0 & 0 & 1 & 1 & 1 \\
\hline 0 & 0 & 0 & 1 & {\color{blue} 0} \\
\hline
\end{array}
\]
\label{e:t}
\end{ex}
It is well known that two formulae are logically equivalent if
they define the same propositional function. Consequently 
they must have the same truth table. Our concern is with 
the set of propositional functions defined by bracketed implications. 
The following uniqueness lemma shows that it
suffices to work with the formulae. 

\begin{lemma}
Two bracketed implications are logically equivalent if and only if
they are equal.
\end{lemma}

\begin{pf}
We show how to recover the bracketing
from the propositional function defined by such a formula. Our proof is
by induction on $n$, the result is trivial for $n\le2$. Suppose that
the proposition function defined by a formula on $t$ distinct variables 
$p_1,\ldots,p_t$, where $1\leq t<n$, recovers the bracketing.

Let $\phi$ be a bracketed implication. Let valuations $\nu_i$
and $\nu_{i,j}$ be defined by

\[\nu_i(p_j) =\cases{0 & if $j=i$,\cr\cr
1 & otherwise.\cr}\;\;\;\;\;\;\;\nu_{i,j}(p_k) =\cases{0 & if $k=i$ or $k=j$, \cr\cr
1 & otherwise.\cr}\]
Now it is straightforward to check that $\nu_n(\phi)=0$, while $\nu_i(\phi)=1$
for $i\ne n$.

Suppose that $\phi$ has the form $\psi\to\chi$, where $\psi$ and $\chi$ are
bracketed implications  involving $p_1,\ldots,p_r$ and $p_{r+1},\ldots,p_n$
respectively. Then, for $i\le r$, we have $\nu_{i,n}(\chi)=0$, while
$\nu_{i,n}(\psi)=1$ if $i<r$, $\nu_{r,n}(\psi)=0$. We conclude that
$\nu_{i,n}(\phi)=0$ if $i<r$ while $\nu_{r,n}(\phi)=1$. Hence we can determine
the value of $r$. By the induction hypothesis, the bracketings of $\psi$ and
$\chi$ are determined by the propositional function, and hence the bracketing
of $\phi$ is determined.\qed
\end{pf}

We could also consider \emph{permuted bracketed implications}, which are formulae
obtained from $p_1\to p_2\to\cdots\to p_n$ by permuting the propositions 
and then inserting brackets, where $p_1,\ldots,p_n$ are distinct propositions. 
More precisely: these are well-formed bracketings 
of $p_{i_1}\to p_{i_2}\to\cdots\to p_{i_n}$, where
$(i_1,\ldots,i_n)$ is a permutation of $(1,\ldots,n)$. Here the situation
is less satisfactory; we can count formulae, but the analogue of our
uniqueness lemma does not hold (for example, $p_1\to(p_2\to p_3)$ and
$p_2\to(p_1\to p_3)$ define the same propositional function), and we do not
know how to count propositional functions represented by permuted bracketed
implications, or the rows with value ``false'' in the corresponding truth
tables.

\section{The number of false rows}

It is well known that the number of bracketings of a product of $n$ terms
is the Catalan number
\[C_n= \frac{1}{n}{2n-2\choose n-1}, \textit{ with } C_0=0\]
whose generating function is
\[\sum_{n\ge1}C_nx^n = (1-\sqrt{1-4x})/2\]
(see~\cite[page 61]{PJCC}).
Then $C_n$ is the number of bracketed implications in $n$ propositional
variables, and by the uniqueness lemma of the preceding section, it is also
the number of propositional functions or truth tables defined by such formulae.

\begin{prop} 
Let $f_n$ be number of rows with the value ``false'' in the truth tables of
all bracketed implications with $n$ distinct variables. Then 
\[f_n = \sum_{i=1}^{n-1} (2^iC_i -f_i) f_{n-i}, \; \textit{ with }\; f_0=0,\; f_1=1. 
\]
\label{p:rec}
\end{prop}

\begin{pf}
A row with the value false comes from an expression $\psi\to \chi$ where $\nu(\psi)=1$
and $\nu(\chi)=0$. If $\psi$ contains $i$ variables, then $\chi$ contains
$n-i$, and the number of choices is given by the summand in the proposition.\qed
\end{pf}
\begin{ex}
\[
f_1=1,\; f_2= (2^1C_1-f_1)f_1=1,
\]
and 
\[
 f_3 = (2^1C_1-f_1)f_2+(2^2C_2-f_2)f_1 = 1+3=4
\]
which coincides with the result we had from Example~\ref{e:t}.
\end{ex}
Using this Proposition, it is straightforward to calculate the values of $f_n$
for small $n$. The first $22$ values are
\begin{eqnarray*}
\{f_n\}_{n\geq 1} &=& 1, 1, 4, 19, 104, 614, 3816, 24595, 162896, 1101922, 7580904, \\
&& 52878654, 373100272, 2658188524, 19096607120, 138182654595, \\
&& 1006202473888, 7367648586954, 54214472633064, \\
&& 400698865376842, 2973344993337520, 22142778865313364, \ldots
\end{eqnarray*}

Let $g_n$ be the total number of rows in all truth tables 
for bracketed implications with $n$ variables. It is clear that $g_n=2^nC_n$, with $g_0=0$.
Let $F(x)$ and $G(x)$ be the generating functions for $f_n$, and $g_n$, respectively. 
That is, $F(x)= \sum_{n\geq 1} f_nx^n$, and $G(x)=\sum_{n\geq 1} g_nx^n$ .
Then Proposition~\ref{p:rec} gives
\begin{equation}\label{eq:1}
F(x) = x + F(x)( G(x) - F(x))
\end{equation}
where $G(x)$ can be obtained from the generating function of $C_n$ by replacing
$x$ by $2x$: that is,
\begin{equation}\label{eq:2}
G(x) = (1- \sqrt{1-8x})/2.
\end{equation}
Substituting the equation~(\ref{eq:2}) into the equation~(\ref{eq:1}) 
gives the following quadratic equation:
\begin{equation}\label{eq:3}
2F(x)^2 + F(x) \left( 1 + \sqrt{1-8x}\right) -2x =0.
\end{equation}
Solving equation~(\ref{eq:3}) gives the following proposition:
\begin{prop}
The generating function for the sequence $\{f_n\}_{n\geq 1}$ is given by
\[
F(x) = \frac{-1-\sqrt{1-8x} + \sqrt{2+2\sqrt{1-8x}+8x}}{4}.
\]
\label{p:gf}
\end{prop}
(As with the Catalan numbers, the choice of sign in the square root is made to
ensure that $F(0)=0$.)
With the help of Maple we can obtain the first $22$ terms of the above series, 
and hence give the first $22$ values of $f_n$; these agree with the values
found from the recurrence relation.

\section{Asymptotic analysis}

In this section we want to get an asymptotic formula for the coefficients of
the generating function $F(x)$ from Proposition~\ref{p:gf}. 
We use the following result~\cite[page 389]{FAC}:

\begin{prop}
Let $a_n$ be a sequence whose terms are positive for sufficiently large $n$.
Suppose that $A(x)=\sum_{n\geq 0} a_nx^n$ converges for some value of $x>0$.
Let $f(x)= (-\ln(1-x/r))^b(1-x/r)^c$, where $c$ is not a
positive integer, and we do not have $b=0$ and $c=0$. 
Suppose that $A(x)$ and $f(x)$ each have a singularity at $x=r$ 
and that $A(x)$ has no singularities in the interval $[-r,r)$. 
Suppose further that $\lim_{x\to r} \frac{A(x)}{f(x)}$ exists and has nonzero value $\gamma$. Then 
\[a_n \sim \cases{
\gamma{n-c-1\choose n} (\ln{n})^br^{-n}, & if $c\not=0$,\cr\cr
\frac{\gamma b(\ln{n})^{b-1}}{n}, & if $c=0$.\cr}\]
\label{p:b}
\end{prop}

\begin{no}\label{n:b}
We also have
\[
{n-c-1\choose n}\sim \frac{n^{-c-1}}{\Gamma(-c)},
\]
where the standard gamma-function
\[\Gamma(x)=\int_0^\infty t^{x-1}\mathrm{e}^{-t}\,\mathrm{d}t\]
satisfies $\Gamma(x+1)=x\Gamma(x)$ and $\Gamma(1/2)=\sqrt{\pi}$. It
follows that $\Gamma(-1/2)=-\sqrt{\pi}/2$.
\end{no}
\medskip

Recall that $G(x)=(1-\sqrt{1-8x})/2$, therefore 
\[F(x)= \frac{(G(x)-1)+\sqrt{(1-G(x))^2+4x}}{2}.\]

Before studying $F(x)$, we first study $G(x)$. This $G(x)$ could easily be
studied by using the explicit formula for its coefficients, which is
$2^n{2n-2\choose n-1}/n$. But our aim is to understand how to handle the
square root singularity. A square root singularity occurs while attempting
 to raise zero to a power which is not a positive integer. 
Clearly the square root, $\sqrt{1-8x}$, has a singularity at $1/8$.
Therefore by Proposition~\ref{p:b}, $r=1/8$. We have
$G(1/8)=1/2$, so we would not be able to divide $G(x)$ by a
suitable $f(x)$ as required in Proposition~\ref{p:b}. To create a function
which vanishes at $\frac{1}{8}$, we simply look at
$A(x)=G(x)-1/2$ instead. That is, let
\[f(x)=(1-x/r)^{1/2} =  (1-8x)^{1/2}.\]
Then
\[\gamma=\lim_{x\rightarrow1/8}\frac{A(x)}{\sqrt{1-8x}} = -\frac{1}{2}.\]

Now by using Proposition~\ref{p:b} and Note~\ref{n:b},
\[
g_n \sim -\frac{1}{2}{n-\frac{3}{2}\choose n} \bigg(\frac{1}{8}\bigg)^{-n}  \sim -\frac{1}{2}\,\frac{8^n n^{-3/2}}{\Gamma (-1/2)} = \frac{2^{3n-2}}{\sqrt{\pi n^3}}.
\]

We are now ready to tackle $F(x)$, and state the main theorem of the paper.

\begin{theorem}
Let $f_n$ be number of rows with the value false in the truth tables of all the
bracketed implications with $n$ distinct variables. Then
\[f_n \sim \left(\frac{3-\sqrt{3}}{6}\right)\frac{2^{3n-2}}{\sqrt{\pi n^3}}.\]
\end{theorem}

\begin{pf} We have
\[F(x) = \frac{-1-\sqrt{1-8x} + \sqrt{2+2\sqrt{1-8x}+8x}}{4}.\]

We find that $r=\frac{1}{8}$, and $f(x)=\sqrt{1-8x}$. 
Since $F(1/8)=(-1+\sqrt{3})/4 \not= 0$, we need a function 
which vanishes at $F(1/8)$, thus we let $A(x)=F(x)-F(1/8)$.

\[\lim_{x\to1/8} \frac{A(x)}{f(x)} 
= \lim_{x\to 1/8} \frac{-\sqrt{1-8x}+\sqrt{2+2\sqrt{1-8x} + 8x} 
-\sqrt{3}}{4\sqrt{1-8x}} .\]

Let $v=\sqrt{1-8x}$. Then
\begin{eqnarray*}
\gamma &=& \lim_{v\to 0} \frac{-v +\sqrt{(1+v)(3-v)} - \sqrt{3}}{4v} \\
&=& \lim_{v\to 0} \frac{-v +\sqrt{3+2v-v^2} - \sqrt{3}}{4v} \\
&=& \lim_{v\to 0} \frac{-1 +(1-v)(3+2v-v^2)^{-1/2}}{4} \\
&=& -\frac{3-\sqrt{3}}{12},
\end{eqnarray*}
where we have used l'H\^opital's Rule in the penultimate line.

Finally,
\[f_n \sim -\frac{3-\sqrt{3}}{12}{n-\frac{3}{2}\choose n}
\left(\frac{1}{8}\right)^{-n} \sim \left(\frac{3-\sqrt{3}}{6}\right)
\frac{2^{3n-2}}{\sqrt{\pi n^3}},\]
and the proof is finished. \qed
\end{pf}

The importance of the constant $(3-\sqrt{3})/6 = 0.2113248654$ 
lies in the following fact:

\begin{cor}
Let $g_n$ be the total number of rows in all truth tables for bracketed
implications with $n$ variables, and $f_n$ the number of rows with the
value ``false''. Then $\lim_{n\to\infty}f_n/g_n=(3-\sqrt{3})/6$.
\end{cor}

The table below illustrates the convergence.

\[
\begin{array}{|c|c|c|c|}
\hline n & f_n & g_n & f_n/g_n \\
\hline 1 & 1 & 2 & 0.5 \\
\hline 2 & 1 & 4 & 0.25 \\
\hline 3 & 4 & 16 & 0.25 \\
\hline 4 & 19 & 80 & 0.2375 \\
\hline 5 & 104 & 448 & 0.2321428571 \\
\hline 6 & 614 & 2688 & 0.228422619 \\
\hline 7 & 3816 & 16896 & 0.2258522727 \\
\hline 8 & 24595 & 109824 & 0.2239492279 \\
\hline 9 & 162896 & 732160 & 0.2224868881 \\
\hline 10 & 1101922 & 4978688 & 0.2213277876 \\
\hline
\end{array}\]
For $n=100$ the ratio is $0.2122908650$, and for $n=1000$ it is
$0.2114211279$.

\begin{cor} \label{c:c}
Let $t_n$ be the number of rows with the value ``true'' in the truth tables of all 
bracketed formulae with $n$ distinct variables connected by the binary 
connective of implication. Then 
\[
t_n= g_n -f_n, \textit{ with } t_0 =0,
\]
and for large $n$, 
\[
t_n \sim \Bigg(\frac{3+\sqrt{3}}{6}\Bigg)\frac{2^{3n-2}}{\sqrt{\pi n^3}} .
\]
\end{cor}
Using this Corollary~\ref{c:c}, it is straightforward to calculate the values of $t_n$. 
The table below illustrates this up to $n=10$.
\[
\begin{array}{|l|c|c|c|c|c|c|c|c|c|c|c|}
\hline n&0 & 1& 2 & 3 & 4 & 5 & 6 & 7 & 8 & 9 & 10\\
\hline t_n & 0 & 1 & 3 & 12  & 61 & 344 & 2074 & 13080 & 85229 & 569264 & 3876766\\
\hline
\end{array}
\]

\end{document}